\documentclass[11pt]{article}

\usepackage{amscd}
\usepackage{amsmath}
\usepackage{amsfonts}
\usepackage{pictex}
\usepackage{theorem}
\usepackage{mathrsfs}

\makeindex

{\theorembodyfont{\slshape}
\newtheorem{theorem}{Theorem}
\newtheorem{proposition}{Proposition}
\newtheorem{lemma}{Lemma}
\newtheorem{corollary}{Corollary}

}
{\theorembodyfont{\rmfamily}
\newtheorem{definition}{Definition}

\newtheorem{remark}{Remark}

}
\def\proof{{\noindent\sc Proof. \quad}}

\def\eproof{{\mbox{}\hfill\qed}\medskip}
\newcommand\qed{{\unskip\nobreak\hfil\penalty50\hskip2em\vadjust{}
\nobreak\hfil$\Box$\parfillskip=0pt\finalhyphendemerits=0\par}}
\def\psid{\psi_{\delta}}
\def\psip{\dot{\psi}}
\begin{document}

\bibliographystyle{plain}

\makeatletter

\renewcommand{\bar}{\overline}

\newcommand{\x}{\times}
\newcommand{\<}{\langle}
\renewcommand{\>}{\rangle}
\newcommand{\into}{\hookrightarrow}

\renewcommand{\a}{\alpha}
\renewcommand{\b}{\beta}
\renewcommand{\d}{\delta}
\newcommand{\D}{\Delta}
\newcommand{\e}{\varepsilon}
\newcommand{\g}{\gamma}
\newcommand{\G}{\Gamma}
\renewcommand{\l}{\lambda}
\renewcommand{\L}{\Lambda}
\newcommand{\n}{\nabla}
\newcommand{\var}{\varphi}
\newcommand{\s}{\sigma}
\newcommand{\Sig}{\Sigma}
\renewcommand{\t}{\tau}
\renewcommand{\O}{\Omega}
\renewcommand{\o}{\omega}
\newcommand{\z}{\zeta}

\newcommand{\p}{\partial}
\renewcommand{\hat}{\widehat}
\renewcommand{\bar}{\overline}
\renewcommand{\tilde}{\widetilde}


\def\N{\mathbb N}
\def\R{\mathbb R}
\def\C{\mathbb C}
\def\E{\mathbb E}


\def\diag{{\sf diag}}
\def\Id{{\rm Id}}
\def\supp{{\rm supp}}
\def\Oh{{\cal O}}
\def\bC{{\bf C}}
\def\DH{{\Delta_\HH}}
\def\DW{{\Delta_W}}
\def\CC{{\mathscr C}}
\def\CH{C_{\HH}}
\def\DD{{\cal D}}
\def\BB{{\cal B}}
\def\AA{{\cal A}}
\def\SS{{\cal S}}
\def\VCD{{\rm VCD}\,}
\def\Range{{\sf Range}}
\def\Var{{\sf Var}}
\def\comb{{\Huge\sf C}}
\def\Hd{{\cal H}_d}
\def\HK{{\cal H}_K}
\def\Fd{{\cal F}_d}
\def\Pd{{\cal P}_d}
\def\PP{{\cal P}}
\def\LL{{\cal L}}
\def\NN{{\cal N}}
\def\Esp{\mathop{\bf E}}
\def\argmin{\mathop{\rm argmin}}
\def\bfE{{\bf E}}
\def\bfe{{\bf e}}
\def\bM{{\bf M}}
\def\bB{{\bf B}}
\def\Ez{{\bf E}_\bz}
\def\Er{{\rm E}}
\def\FF{{\cal F}}
\def\cG{{\cal G}}
\def\EE{{\cal E}}
\def\EH{{\cal E}_{\HH}}
\def\EHz{{\cal E}_{\HH,\bz}}
\def\tE{{\cal E}_\rho}
\def\ID{{\cal I}_\rho}
\def\HH{{\cal H}}
\def\Prob{\mathop{\sf Prob}}
\def\bz{{\bf z}}
\def\bx{{\boldsymbol x}}
\def\bm{{\boldsymbol m}}
\def\bv{{\boldsymbol v}}
\def\ba{{\boldsymbol a}}
\def\bb{{\boldsymbol b}}
\def\bH{{\boldsymbol H}}
\def\bPhi{{\boldsymbol \Phi}}
\def\L2{{\mathscr L}^2_\rho(X)}
\def\LL{{\mathscr L}}


\def\JACM{Journal of the ACM}
\def\CACM{Communications of the ACM}
\def\ICALP{International Colloquium on Automata, Languages
            and Programming}
\def\STOC{annual ACM Symp. on the Theory
          of Computing}
\def\FOCS{annual IEEE Symp. on Foundations of Computer Science}
\def\SIAM{SIAM Journal on Computing}
\def\BSMF{Bulletin de la Soci\'et\'e Ma\-th\'e\-ma\-tique de France}
\def\CRAS{C. R. Acad. Sci. Paris}
\def\IPL{Information Processing Letters}
\def\TCS{Theoretical Computer Science}
\def\BAMS{Bulletin of the Amer. Math. Soc.}
\def\TAMS{Transactions of the Amer. Math. Soc.}
\def\PAMS{Proceedings of the Amer. Math. Soc.}
\def\JAMS{Journal of the Amer. Math. Soc.}
\def\LNM{Lect. Notes in Math.}
\def\LNCS{Lect. Notes in Comp. Sci.}
\def\JSL{Journal for Symbolic Logic}
\def\JSC{Journal of Symbolic Computation}
\def\JCSS{J. Comput. System Sci.}
\def\JoC{J. of Complexity}
\def\MP{Math. Program.}
\sloppy

\begin{title}
{\Large {\bf Flocking in noisy environments}}
\end{title}
\author{Felipe Cucker\\
Department of Mathematics\\
City University of Hong Kong\\
83 Tat Chee Avenue, Kowloon\\
HONG KONG\\
e-mail: {\tt macucker@math.cityu.edu.hk}
\and
Ernesto Mordecki\\
Centro de Matem\'{a}tica.\ Facultad de Ciencias.\\
Universidad de la Rep\'{u}blica\\
Igu\'{a} 4225. 11400 Montevideo\\
URUGUAY\\
e-mail:  {\tt mordecki@cmat.edu.uy} 
}

\date{}
\makeatletter
\maketitle
\makeatother

\section{Introduction}\label{sec:intro}

The problem of reaching a consensus in a group of autonomous 
agents has been the object of study in a number of 
situations ranging from 
linguistics~\cite{CSZ:04,WangEtAl01,Niyogi:06} 
to distributed computing~\cite{Tsitsi:84,TBA:86} and from 
physics~\cite{VCBC:95} to animal 
behavior~\cite{CHDB:06,FGLO:99,TJP:03ab,TBL:06}.  

Research on the latter attempts to explain, by appropriately 
modeling it, the observed behavior of a group of animals, say 
a flock of birds, whose velocities converge to a common one. 
An influential 
model for this behavior has been postulated in~\cite{VCBC:95} 
by Vicsek and collaborators and studied in~\cite{JLMo:03} where 
convergence is shown under some conditions on the sequence of 
states of the flock. Also, the work of Buhl et al.~\cite{Buhl:06} 
showed that Vicsek's model accurately describes the convergence 
to order in swarms of locusts. A different model for the same 
phenomenon was proposed in~\cite{CS06} and extended 
in~\cite{CS07} (to situations other than flocking) and 
in~\cite{Shen:07} (to flocking situations where a hierarchical 
leadership structure is present). A main feature of these 
papers is 
that, in contrast with the results in~\cite{JLMo:03,VCBC:95}, 
convergence results are established conditioned to the {\em initial} 
state of the flock only. 

The model in~\cite{CS06} postulates the following behavior: 
every bird adjusts 
its velocity by adding to it a weighted average of the differences 
of its velocity with those of the other birds. That is, at time 
$t$, and for bird $i$,
\begin{equation}\label{eq:model}
   \bv_i(t+h)=\bv_i(t)+h\sum_{j=1}^ka_{ij}(\bv_j(t)-\bv_i(t))
\end{equation}
where the weights $\{a_{ij}\}$ quantify the way the birds 
influence 
each other, and $h>0$ is the time step. Vicsek's model is more 
general in the sense that allows the obtained $\bv_i(t+h)$ 
to be perturbed  by some centered noise. More precisely, 
it replaces 
(\ref{eq:model}) by 
\begin{equation}\label{eq:mod_noise}
   \bv_i(t+h)=\bv_i(t)+h\sum_{j=1}^ka_{ij}(\bv_j(t)-\bv_i(t))
   +h \bH_i(t)
\end{equation}
where $\bH_i(t)\in\E$ is some centered random variable 
modeling the noise. 
Here $\E$ denotes 3-dimensional Euclidean space. 
 
It is reasonable to assume that the influence weights $a_{ij}$ are 
a function of the distance between birds. This is the case 
in~\cite{CS06,VCBC:95} and the major difference between these 
models is in the choice of the $a_{ij}$. In this paper we will 
follow~\cite{CS06} but slightly depart from it and take 
the {\em adjacency matrix} $A_x$ to have entries 
\begin{equation}\label{eq:adj}
  a_{ij}=\frac{K}{(1+\|\bx_i-\bx_j\|)^\a}
\end{equation} 
for some fixed $K>0$ and $\a\geq 0$.

We can write the set of equalities (\ref{eq:model}) in a more 
concise form. Let $D_x$ be the $k\x k$ diagonal matrix 
whose $i$th diagonal entry is $d_i=\sum_{j\leq k} a_{ij}$ and 
$L_x=D_x-A_x$. Then (cf.~\cite{CS06})
$$
   \bv(t+h)-\bv(t)=-hL_x \bv(t)+h\bH(t)
$$ 
where $\bH=(\bH_1,\ldots,\bH_k)$.  

Note that the matrix notation $L_x \bv(t)$ does not have the usual 
meaning of a $k\times k$ matrix acting on $\R^k$. Instead, the 
matrix $L_x$ is acting on $\E^k$ by mapping 
$(\bv_1,\ldots,\bv_k)$ to 
$((L_x)_{i1}\bv_1+\ldots+(L_x)_{ik}\bv_k)_{i\leq k}$. 

Adding a natural equation for the change of positions 
we obtain the discrete dynamical system
\begin{align}\label{eq:our}
  \bx(t+h)\;&=\;\bx(t)+h \bv(t)\tag{D}\\
  \bv(t+h)\;&=\;\left(\Id- hL_x\right)\bv(t) +h\bH(t).\notag
\end{align}

We also consider evolution for continuous time. 
The corresponding model is obtained by taking limits 
for $h\to0$ and can be given by the 
system of differential equations
\begin{align}\label{eq:dif_cont}
\bx'\;&=\; \bv\tag{C}\\ 
\bv'\;&=\; -L_x \bv + \bH.\notag
\end{align}

The main result in~\cite{CS06} shows, 
in the unperturbed case, for both discrete and 
continuous time, convergence to the alignment of the velocities. 
More precisely, convergence to a common velocity when the 
initial positions and velocities of the flock are not too 
dissimilar (for $\a\geq1$, otherwise, convergence holds 
unconditionally). For systems~\eqref{eq:our} 
and~\eqref{eq:dif_cont}, due to the presence of noise, 
we can not expect convergence to a common 
velocity. Once the velocities $\{\bv_1,\ldots,\bv_k\}$ are 
similar enough compared with the noise  
the latter will, with positive probability, outdo the 
contractive character of the system. 
Perfect alignment as in~\cite{CS06} should therefore be replaced  
by ``nearly-alignment.'' A formal measure of similarity  
(and with it a definition of nearly-alignment) will be given 
soon in Section~\ref{sec:prelim}. A description of the forms of 
noise we consider in this paper will be given 
in \S\ref{sec:stat_disc} and \S\ref{sec:stat_cont}.
We nevertheless state now an informal version 
of our main results (see Theorems~\ref{th:disc_absolute} 
and~\ref{th:main_cont} for precise statements). 
\medskip

\noindent
{\bf Main result}\quad
{\sl Assume that at time 0 the positions and velocities of the flock are not both too dissimilar  
(for $\a\geq 1$, otherwise, no assumption is needed) and that 
the time step $h$ is small enough (in case of discrete time). 
Then nearly-alignment is (quickly) reached with a certain probability
and we exhibit a lower bound for this probability in terms of 
the initial similarity of positions and velocities, the variance 
of $\bH$, and the parameters $k,K$ and $\a$. }
\medskip

The proof of our main result closely follows the proofs 
in~\cite{CS06}. Some changes had to be made to make room 
for the noise and in doing so we did a few simplifications 
as well.  

\section{Some preliminaries}\label{sec:prelim}

\subsection{Laplacians}

Given a nonnegative, symmetric, $k\x k$ matrix $A$ the 
{\em Laplacian} $L$ of $A$ is defined to be
$$
   L=D-A
$$ 
where $D=\diag(d_1,\ldots,d_k)$ and $d_\ell=\sum_{j=1}^ka_{\ell j}$. 
Some features of $L$ are immediate. It is symmetric 
and it does not depend on the diagonal entries of $A$. 

The matrix $L_x$ in (\ref{eq:our}) and (\ref{eq:dif_cont}) is 
thus the Laplacian of $A_x$. It satisfies that 
for all $u\in\E$, $L_x(u,\ldots,u)=0$. In addition, 
it is positive semidefinite. 

The smallest eigenvalue of $L_x$ is zero. 
Its second eigenvalue is called the {\em Fiedler 
number} of $A_x$. We will denote it by $\phi_x$. 

\begin{proposition}{\bf (\cite[Proposition~1]{CS07})}\label{prop:phi}
Let $A$ be a $k\times k$ nonnegative, symmetric matrix, 
$L=D-A$ its Laplacian, $\phi$ its Fiedler number, 
and $\mu=\min_{i\neq j} a_{ij}$. 
Then $\displaystyle \phi\geq k\mu$. 
\eproof
\end{proposition}

\subsection{Similarity and nearly-alignment}

The inner product on $\E$ naturally induces an inner product 
on $\E^k$. 
Let $\Delta$ be the diagonal of $\E^k$, i.e.,
$$
  \Delta=\{(u,u,\ldots,u)\mid u\in\E\}
$$
and $\Delta^\perp$ be the orthogonal complement of $\Delta$ in 
$\E^k$. Then, every point $\bv\in\E^k$ decomposes in a 
unique way as $\bv=\bv_\Delta+\bv_\perp$ with 
$\bv_\Delta\in\Delta$ and $\bv_\perp\in\Delta^\perp$. 
This decomposition has a simple explicit form. Denote 
by 
$$
  \bm=\frac{1}{k}\sum_{i=1}^k \bv_i
$$
the mean of the $\bv_i$. Then $\bv_{\Delta}=(\bm,\ldots,\bm)$ 
and $\bv_{\perp}=(\bv_1-\bm,\ldots,\bv_k-\bm)$. This follows 
immediately from the equality
\begin{eqnarray*}
  \langle \bv_{\Delta},\bv_{\perp}\rangle 
  &=&\sum_{i=1}^k \langle \bm, (\bv_i-\bm)\rangle 
  =\left\langle \bm, \sum_{i=1}^k \bv_i-\bm\right\rangle \\
  &=&\left\langle \bm, 
    \left(\sum_{i=1}^k \bv_i\right)-k\bm\right\rangle 
  =\langle \bm, 0\rangle =0.
\end{eqnarray*}
We can look at the evolution of the velocities $\bv_i(t)$ 
decomposing into the evolution of their mean $\bm(t)$  
and that of the distances to that mean 
$\bv_{\perp}=(\bv_1-\bm,\ldots,\bv_k-\bm)$ and a key observation 
at this stage is the fact that convergence to a common velocity 
(or nearly-alignment) is a feature of the second evolution only. 
More precisely, the condition 
``the velocities $\bv_i(t)$ tend to alignment 
is equivalent to the condition ``$\bv_\perp(t)\to0$''. 
We are thus interested on the projection 
$(\bx_{\perp}(t),\bv_{\perp}(t))$ over 
$\Delta^\perp\x\Delta^\perp$ of 
the solutions $(\bx(t),\bv(t))$ 
of the system (\ref{eq:our}) (or (\ref{eq:dif_cont})). 
It is easy to show (see~\cite{CS06}) that 
these projections are the solutions of the restriction 
of (\ref{eq:our}) (resp.~(\ref{eq:dif_cont})) 
to $\Delta^\perp\x\Delta^\perp$. 

More precisely, they are the solutions of 
\begin{align*}
  \bx(t+h)_\perp\;&=\;\bx(t)_\perp+ h \bv(t)_\perp\\
  \bv(t+h)_\perp\;&=\;\left(\Id- hL_{x_\perp}\right)\bv(t)_\perp 
  +h\bH_\perp.
\end{align*}
Hence, in what follows, we will consider positions in 
$$
    X:=\E^k/\Delta\simeq \Delta^\perp
$$
and velocities in
$$
    V:=\E^k/\Delta\simeq \Delta^\perp.
$$
For $\bx,\bv\in \E^k$ we will denote $x=\bx_\perp$ and 
$v=\bv_\perp$. Finally, we will denote $H=\bH_\perp$.  

It is natural now to take the norm $\|\bx_\perp\|$ 
of the projection $\bx_\perp$ as the {\em dissimilarity} of 
$\bx$ and similarly for $\|\bv_\perp\|$. In the case 
of $\bx$ we may call this measure the {\em dispersion} of the 
flock. It relates with its ``diameter.'' 

\begin{lemma}\label{lem:gamma}
For all $\bx\in\E^k$, 
$\max_{i\neq j}\|\bx_i-\bx_j\|\leq \sqrt{2}\|\bx_\perp\|$.
\end{lemma}

\proof
Write $\bx=\bx_\Delta +\bx_\perp =
 (\tilde u,\ldots,\tilde u) + ((\bx_\perp)_1,\ldots,(\bx_\perp)_k)$. 
Then, for all $i\neq j$, 
$\bx_i-\bx_j=(\bx_\perp)_i-(\bx_\perp)_j$ and 
\begin{equation}\tag*{\qed}
 \|\bx_i-\bx_j\|_{\E}=\|(\bx_\perp)_i-(\bx_\perp)_j\|_{\E} 
 \leq \|(\bx_\perp)_i\|_{\E}+\|(\bx_\perp)_j\|_{\E} 
 \leq \sqrt{2}\|\bx_\perp\|_{\E^k}.
\end{equation}

The notion of similarity leads to the following definition.

\begin{definition}
Let $\nu>0$. We say that the flock $\{1,\ldots,k\}$ is 
{\em $\nu$-nearly-aligned} (or simply {\em nearly-aligned}) when 
$\|\bv_\perp\|\leq\nu$. 
\end{definition} 

\subsection{A few functions of the initial state}

The initial state of the flock is characterized by the 
pair $(\bx(0),\bv(0))$. For convergence to alignment 
(or to nearly-alignment) to hold one needs to require 
that the dissimilarities of these two vectors are not both 
large. 

We close this section with a few quantities related to these 
initial dissimilarities which will occur when describing the 
conditions ensuring convergence. These are 
$$
  \ba=\frac{2\sqrt{2}}{kK}\|v(0)\|, 
  \qquad
  \bb=1+\sqrt{2}\|x(0)\|, 
$$
$$
  U_0=\left\{\begin{array}{ll}
  \max\left\{
  \left(2\ba\right)^{\frac{1}{1-\a}}, 
  2\bb\right\}
  & \mbox{ if $\a<1$}\\ [10pt]
  \frac{\bb}{1-\ba} 
  & \mbox{ if $\a=1$}\\ [10pt]
  \frac{\a}{\a-1}\bb    
  & \mbox{ if $\a>1$,}
     \end{array}\right.
$$
$$
  B_0=\frac{U_0-1}{\sqrt{2}},\qquad
  \mbox{and}\qquad  
 \HH_0=\frac{2^{-\a-1}kK}{U_0^\a}.
$$

\section{Discrete time}\label{sec:disc_time}

Assume the initial state for \eqref{eq:our} 
is at time 0. Then the sequence of 
states is $\{x(th),v(th)\}_{t\in\N}$. To simplify notation 
we will denote $x(th)$ simply by $x[t]$ and similarly for 
$v$. 

\subsection{Statement of the result}\label{sec:stat_disc}

Recall, the random noise in (\ref{eq:our}) has the form 
$\bH=(\bH_1,\ldots,\bH_k)$ and we have 
$\bH=\bH_\Delta+ H$. Note that the component $\bH_{\Delta}$ 
of $\bH$ corresponds to the perturbation of the 
common velocity $\bv_{\Delta}$ within $\Delta$ and is 
therefore of no consequence regarding convergence to 
alignment or nearly-alignment.  

In what follows we assume that, for all $i\in\{1,\ldots,k\}$, 
and for all $t\in\N$,
$$
   \bH_i(t)=\left(e_i^{(1)}(t),e_i^{(2)}(t),e_i^{(3)}(t)\right)
$$
where the $e_i^{(\ell)}(t)$ are one dimensional random variables,
the coordinates of the perturbation.

We consider two possible laws for the distribution of $\bH$: 
$$
  {\bf Uniform:}\qquad \bH \simeq U_{3k}(0,r)
$$
where, for some $r>0$, 
$U_{3k}(0,r)$ is the uniform distribution in 
$B(0,r)\subset\R^{3k}$, and 
$$
  {\bf Gaussian:}\qquad \bH\simeq N(0,\sigma^2\Id_{3k}),
$$
a $3k$-dimensional 
centered Gaussian distribution with covariance matrix $\sigma^2\Id_{3k}$.
As a consequence,  
in the Gaussian case, the random 
variables $e_i^{(\ell)}$ are independent. 

Our main result for discrete time is the following.

\begin{theorem}\label{th:disc_absolute}
Consider the system \eqref{eq:our} with adjacency 
matrix given by~\eqref{eq:adj}.  
Assume that $h$ satisfies
$$
  h<\min\left\{\frac{1}{2(k-1)\sqrt{k}K},
  \frac{1}{2\sqrt{2}\|v(0)\|}
  \left(\frac{kK}{2\HH_0}\right)^{1/\a}\right\}.
$$
Assume also that one of the three following hypothesis 
holds:
\begin{description}
\item[(i)]
$\a<1$, 
\item[(ii)]
$\a=1$, and $\|v(0)\|< \frac{kK}{2\sqrt{2}}$, 
\item[(iii)]
$\a>1$, and 
$$
   \left(\frac1{\a\ba}\right)^{\frac1{\a-1}}
   \frac{\a-1}{\a} > \bb +2kKh\ba.
$$
\end{description}
Then $\nu$-nearly-alignment for some $\nu<\|v(0)\|$
occurs in a number of iterations bounded by 
$$
  T_0:=\frac{2U_0^\a}{hkK}
  \ln\left(\frac{\|v(0)\|}{\nu}\right)
$$
with probability at least
$$
  \left(\frac{\HH_0\nu}{r}\right)^{3kT_0}
$$
in the uniform case (1 if ${r}\leq\HH_0\nu$), and with 
probability at least
\begin{equation*}\label{eq:gaussian_bound}
  \left(\int_0^{\sqrt{\HH_0\nu/(2{\sigma})}}
  \frac{t^{\frac{3k-5}{2}}}
 {\Gamma\left(\frac{3k-3}{2}\right)}e^{-t}dt\right)^{T_0}
\end{equation*}
in the Gaussian case. 
\end{theorem}

\begin{remark}
For each of the cases (i), (ii), and (iii) we can replace 
$U_0$ and $\HH_0$ by their respective values. In case (iii) 
and with uniform noise, for instance, this yields
$$
  T_0=\frac{2}{hkK}\ln\left(\frac{\|v(0)\|}{\nu}\right)
      \left((1+\sqrt{2}\|x(0)\|)\frac{\a}{\a-1}\right)^\a
$$
and
$$
  \frac{\HH_0\nu}{{r}}
   = \frac{2^{-\a-1} kK}{{r}}
       \left(\frac{\a-1}{\a(1+\sqrt{2}\|x(0)\|)}\right)^\a. 
$$
Note that this means that for 
$$
  T_0=\frac{1}{hkK}\ln\left(\frac{\|v(0)\|}{\nu}\right)
      \Oh(\|x(0)\|^\a)
$$
we have
$$
  \Prob\{\mbox{nearly-align in at most $T_0$ iterations}\}
  \ge \left(\frac{kK}{{r}\Oh(\|x(0)\|^\a)}\right)^{3kT_0}.
$$
From these expressions it is easy to read the role of 
the deterministic setting parameters $h,k$ and $K$, the 
probabilistic 
${r}$, the radius $\nu$, and the initial dissimilarities   
$\|x(0)\|,\|v(0)\|$ both in the time required to reach  
nearly-alignment and in the confidence with which this occurs. 
\end{remark}

\begin{remark}
The integral in the bound for the probability in the 
Gaussian case satisfies, 
when $k$ is odd and writing $n=\frac{3k-3}{2}$, the equality
$$
  \int_0^{x}\frac{t^{n-1}}
  {\Gamma\left(n\right)}e^{-t}dt
 =1-e^{-x}\left(\frac{x^{n-1}}{(n-1)!}+\cdots+\frac x{1!}+1\right).
$$
For $x=\sqrt{\HH_0\nu/(2\sigma)}$ and for small $\sigma$ 
this probability bound is equivalent to 
$$
  1-\frac{T_0}{\Gamma(n)}e^{-x}x^{n-1}.
$$
By L'H\^opital's rule, this equivalence holds as well  
when $k$ is even.
\end{remark}

\subsection{Bounded noise}\label{sec:disc_rob}

Fix a solution $(x,v)$ of (\ref{eq:our}). 
At a time $t\in\N$, $x[t]$ and $v[t]$ are elements 
in $X$ and $V$, respectively. In particular, $x[t]$ determines an 
adjacency matrix $A_{x[t]}$. For notational simplicity 
we will denote 
its Laplacian and Fiedler number by 
$L_t$ and $\phi_t$, respectively. 

\begin{lemma}
For all $x\in X$, 
$$
  \|L_x\|\leq 2(k-1)\sqrt{k}K.
$$
In particular, if $h<\frac{1}{2(k-1)\sqrt{k}K}$ then 
$h\|L_x\|\in(0,1]$. 
\end{lemma}

\proof
For all $i,j\leq k$, $a_{ij}\leq K$. 
Therefore,
$$
  \|L_x\|_{\max}=\max_{i\leq k}\sum_{j=1}^k|(L_x)_{ij}| 
   \leq 2(k-1)K.
$$
Now use that 
$\|L_x\|\leq \sqrt{k}\|L_x\|_{\max}$~\cite[Table~6.2]{Higham96} 
to deduce the result. 
\eproof

\begin{proposition}\label{prop:disc1}
Assume that $h<\frac{1}{2(k-1)\sqrt{k}K}$. 
Assume also that, for all $0\leq t<T$, $\|H\|\leq \HH_0\|v[t]\|$. Then, 
for all $t<T$, 
$$
   \|v[t+1]\|\leq (1-h\phi_t+h\HH_0)\|v[t]\|.
$$
In particular, $\|v\|$ is decreasing as a function 
of $t$ for $t<T$.
\end{proposition}

\proof
The linear map $\Id-hL_t$ is self-adjoint and 
its eigenvalues are in the interval 
$(0,1)$. Its largest eigenvalue is 
$1-h\phi_t$. Therefore  
\begin{align}
 \|v[t+1]\|\;&=\;\left\|\left(\Id-hL_t\right)v[t]+hH\right\|
 \leq \|\Id-hL_t\|\|v[t]\|+h\|H\|\notag\\
 &\leq\; (1-h\phi_t)\|v[t]\|+h\HH_0\|v[t]\|
 = (1-h\phi_t+h\HH_0)\|v[t]\|.\tag*{\qed}
\end{align}

\begin{corollary}\label{prop:disc2}
In the hypothesis of Proposition~\ref{prop:disc1}, for all 
$t\in\{0,\ldots,T-1\}$ we have 
\begin{equation}\tag*{\qed}
  \|v[t]\|\leq \|v[0]\|\prod_{i=0}^{t-1}(1-h\phi_i+h\HH_0).
\end{equation}
\end{corollary}

A proof of the following lemma is in~\cite[Lemma~7]{CS03}.

\begin{lemma}\label{cero}
Let $c_1,c_2>0$ and $s>q>0$. Then the equation 
$$
   F(z)=z^s-c_1z^q-c_2=0
$$
has a unique positive zero $z_*$. In addition 
$$
  z_*\leq\max\left\{(2c_1)^{\frac{1}{s-q}},
         (2c_2)^{\frac{1}{s}}\right\}
$$ 
and $F(z)\leq 0$ for $0\leq z\leq z^*$.
\eproof
\end{lemma}

\begin{theorem}\label{th:2}
Let $T\in\N\cup\{+\infty\}$. 
Assume that, for all $0\leq t<T$, $\|H\|\leq \HH_0\|v[t]\|$, 
and that $h$ satisfies
$$
  h<\min\left\{\frac{1}{2(k-1)\sqrt{k}K},
  \frac{1}{2\sqrt{2}\|v[0]\|}
  \left(\frac{kK}{2\HH_0}\right)^{1/\a}\right\}
$$
where $\HH_0$ is as in Theorem~\ref{th:disc_absolute}. 
Assume also that one of the three following hypothesis 
holds:
\begin{description}
\item[(i)]
$\a<1$, 
\item[(ii)]
$\a=1$, and $\|v[0]\|< \frac{kK}{2\sqrt{2}}$, 
\item[(iii)]
$\a>1$, and
$$
   \left(\frac1{\a\ba}\right)^{\frac1{\a-1}}
   \frac{\a-1}{\a} > \bb +2kKh\ba.
$$
\end{description}
Then $1-h\frac{kK}{2U_0^\a}\in (0,1)$, 
for all $0\leq t <T$, $\|x[t]\|\leq B_0$ and  
$$
 \|v[t]\|\leq \|v[0]\|
 \left(1-h\frac{kK}{2U_0^\a}\right)^t.
$$
In particular, when $T=\infty$, $\|v[t]\|\to0$ for $t\to\infty$.  
\end{theorem}

\proof
Let 
$$
   \Upsilon=\left\{t\in\{0,\ldots,T-1\}\mid 
      \left(1+\sqrt{2}\|x(t)\|\right)^\a 
       \leq \frac{kK}{2\HH_0}\right\}. 
$$
Note that in all three cases ((i), (ii), and (iii)) the 
definition of $\HH_0$ implies that $0\in\Upsilon$ and hence, that 
$\Upsilon\neq\emptyset$. Assume that $\Upsilon\neq \{0,\ldots,T-1\}$ 
and let $\hat t=\min \{\{0,\ldots,T-1\}\setminus\Upsilon\}$. 

For $t<T$, let $t^*$ be the point maximizing $\|x\|$ in 
$\{0,1,\ldots,t\}$. Then, by Proposition~\ref{prop:phi} 
and Lemma~\ref{lem:gamma}, for $i\in\{0,1,\ldots,t\}$, 
$$
  \phi_i
  \geq \frac{kK}{(1 +\sqrt{2}\|x[i]\|)^\a} 
  \geq \frac{kK}{(1 +\sqrt{2}\|x[t^*]\|)^\a}.
$$
Moreover, since $t^*\leq t<\hat t$ we have
$$
  \phi_{i}-\HH_0 \geq 
  \frac{kK}{(1 +\sqrt{2}\|x[t^*]\|)^\a}-\HH_0 
  \geq \frac{kK}{2(1 +\sqrt{2}\|x[t^*]\|)^\a} =:R(t^*)
$$
Using Corollary~\ref{prop:disc2} we obtain, for all $\tau\leq t$,  
\begin{eqnarray*}
  \|x[\tau]\|  
  &\leq&\|x[0]\|+\sum_{j=0}^{\tau-1}
     \|x[j+1]-x[j]\|
  \leq \|x[0]\|+
         h \sum_{j=0}^{\tau-1} \|v[j]\|\\
  &\leq& \|x[0]\|+
         h \left(\|v[0]\|+\sum_{j=1}^{\tau-1} \|v[j]\|\right)\\
  &\leq& \|x[0]\|+
        h \left(\|v[0]\|+\sum_{j=1}^{\tau-1} 
          \|v[0]\|\prod_{i=1}^{j}
          \left(1-h\phi_i+h\HH_0\right)\right)\\
  &\leq& \|x[0]\|+
    h \|v[0]\|\sum_{j=0}^{\tau-1} 
    \left(1-hR(t^*)\right)^{j}\\
  &\leq& \|x[0]\|+
   h \frac{1}{hR(t^*)}\|v[0]\|\\
  &=& \|x[0]\|+ 
   \frac{2(1+\sqrt{2}\|x[t^*]\|)^\a}{kK}
   \|v[0]\|.
\end{eqnarray*}
Multiplying by $\sqrt{2}$ and taking $\tau=t^*$, 
the inequality above takes the following equivalent form
$$
  \sqrt{2}\|x[t^*]\|\leq 
  \sqrt{2}\|x[0]\|
   +\frac{2\sqrt{2}(1+\sqrt{2}\|x[t^*]\|)^\a}{kK}\|v[0]\|
$$
or yet
\begin{equation}\label{eq:disc1O}
  (1+\sqrt{2}\|x[t^*]\|)\leq 
  (1+\sqrt{2}\|x[0]\|) 
   +\frac{2\sqrt{2}(1+\sqrt{2}\|x[t^*]\|)^\a}{kK}\|v[0]\|.
\end{equation}
Let $z=1+\sqrt{2}\|x[t^*]\|$. 
Then (\ref{eq:disc1O}) can be rewritten as $F(z)\leq 0$ with
$$
  F(z)=z-\ba z^{\a}-\bb.
$$
\smallskip

\noindent{\bf (i)\quad} Assume $\a<1$. 
By Lemma~\ref{cero}, $F(z)\leq 0$ implies that 
$(1+\sqrt{2}\|x[t^*]\|)\leq U_0$. 
Since $U_0$ is independent of $t$ we deduce that, for all 
$t<\hat t$, 
$$
   \|x[t]\|\leq \frac{U_0-1}{\sqrt{2}}=B_0.
$$
Therefore, for all $t<\hat t$, 
$$
  (1+\sqrt{2}\|x[t]\|)^\a\leq (1+\sqrt{2}\|x[t^*]\|)^\a
  \leq U_0^\a \leq 2^{-\a}\frac{kK}{2\HH_0}
$$
the last by the definition of $\HH_0$. It follows that 
$$
  \|x[\hat t]\| = \|x[\hat t-1]\|+h\|v[\hat t-1]\|
  \leq \|x[\hat t-1]\|+h\|v[0]\|
$$
and therefore
\begin{eqnarray*}
  (1+\sqrt{2}\|x[\hat t]\|) &\leq &
   1+\sqrt{2}\|x[\hat t-1]\|+\sqrt{2}h\|v[0]\|\\
  &\leq& \left(2^{-\a}\frac{kK}{2\HH_0}\right)^{1/\a} +
  \sqrt{2}h\|v[0]\| \leq \left(\frac{kK}{2\HH_0}\right)^{1/\a}
\end{eqnarray*}
the last by our hypothesis on $h$. This is in contradiction with 
the definition of $\hat t$ and shows that no such $\hat t$ exists. 
That is, for all $t<T$, $\|x[t]\|\leq B_0$ and 
$$
   \phi_t-\HH_0\geq F_0
   :=\frac{kK}{2(1+\sqrt{2}B_0)^\a}=\frac{kK}{2U_0^{\a}}.
$$
By Corollary~\ref{prop:disc2}, for $t<T$, 
$$
 \|v[t]\|\leq \|v[0]\|\prod_{i=0}^{t-1}
 \left(1-h\phi_i+h\HH_0\right)
 \leq \left(1-hF_0\right)^t\|v[0]\|.
$$
The convergence results for the case $T=\infty$ now 
readily follow (cf.~\cite[Theorem~3]{CS06}).
\smallskip

\noindent{\bf (ii)\quad} Assume now $\a=1$. Then 
(\ref{eq:disc1O}) takes the form 
$$
  (1+\sqrt{2}\|x[t^*]\|)\left(1- 
   \frac{2\sqrt{2}}{kK}\|v[0]\|\right)
 -\left(1+\sqrt{2}\|x[0]\|\right)\leq 0.
$$
By hypothesis, $\|v[0]\|< \frac{kK}{2\sqrt{2}}$. This 
implies that  
$$
  \|x[t^*]\|\leq 
  kK\frac{1+\sqrt{2}\|x[0]\|}
  {\sqrt{2}kK-4\|v[0]\|}-1=B_0. 
$$
We conclude that, for all $t<\hat t$, 
$$
  (1+\sqrt{2}\|x[t]\|)^\a\leq 1+\sqrt{2}\|x[t^*]\|
  \leq kK\frac{1+\sqrt{2}\|x[0]\|}
  {kK-2\sqrt{2}\|v[0]\|} \leq \frac{kK}{4\HH_0}
$$
by the definition of $\HH_0$.
We now proceed as in case~(i). 
\smallskip

\noindent{\bf (iii)\quad} Assume finally $\a>1$. 
The derivative $F'(z)=1-\a\ba z^{\a-1}$ has a unique zero at 
$z_*=\left(\frac1{\a\ba}\right)^{\frac1{\a-1}}$ and 
\begin{eqnarray*}
 F(z_*) &=&\left(\frac1{\a\ba}\right)^{\frac1{\a-1}}
        -\ba\left(\frac1{\a\ba}\right)^{\frac{\a}{\a-1}}-\bb\\
     &=& \left(\frac1{\a}\right)^{\frac1{\a-1}}
         \left(\frac1{\ba}\right)^{\frac1{\a-1}}
        -\left(\frac{1}{\a}\right)^{\frac{\a}{\a-1}}
         \left(\frac1{\ba}\right)^{\frac{1}{\a-1}}-\bb\\
     &=& \left(\frac1{\ba}\right)^{\frac1{\a-1}}
         \left(\frac1{\a}\right)^{\frac1{\a-1}}
         \frac{\a-1}{\a}-\bb\\
     &>& 0
\end{eqnarray*}
the last by our hypothesis. Since $F(0)=-\bb<0$, 
$F''(z)=\a(\a-1)\ba z^{\a-2}>0$ for all $z>0$, and $F(z)\to-\infty$ 
when $z\to\infty$, we deduce that the shape of $F$ is as follows:

\begin{center}

\font\thinlinefont=cmr5
\begingroup\makeatletter\ifx\SetFigFont\undefined%
\gdef\SetFigFont#1#2#3#4#5{%
 \reset@font\fontsize{#1}{#2pt}%
 \fontfamily{#3}\fontseries{#4}\fontshape{#5}%
 \selectfont}%
\fi\endgroup%
\mbox{\beginpicture
\setcoordinatesystem units <0.40000cm,0.40000cm>
\unitlength=0.40000cm
\linethickness=1pt
\setplotsymbol ({\makebox(0,0)[l]{\tencirc\symbol{'160}}})
\setshadesymbol ({\thinlinefont .})
\setlinear
%
%
\linethickness= 0.500pt
\setplotsymbol ({\thinlinefont .})
{\putrule from 10.636 16.351 to 10.636 15.875
}%
%
%
\linethickness= 0.500pt
\setplotsymbol ({\thinlinefont .})
{\putrule from  3.175 16.192 to 23.178 16.192
%
%
\plot 22.924 16.129 23.178 16.192 22.924 16.256 /
}%
%
%
\linethickness= 0.500pt
\setplotsymbol ({\thinlinefont .})
{\putrule from  6.668 16.351 to  6.668 15.875
}%
%
%
\linethickness= 0.500pt
\setplotsymbol ({\thinlinefont .})
{\putrule from  4.128  7.620 to  4.128 23.495
%
%
\plot  4.191 23.241  4.128 23.495  4.064 23.241 /
}%
%
%
\put{\small $z_{0}$} [lB] at  8.028 15.334
\put{\small $z_*$} [lB] at  10.368 15.334
\put{\small $(0,-\bb)$} [lB] at  1.168 11.834
\put{\small $\scriptstyle\bullet$} [lB] at  3.968 11.884
\put{\small $(z_*,F(z_*))$} [lB] at  9.168 19.234
\put{\small $\scriptstyle\bullet$} [lB] at  10.458 18.354
\put{\small $z_\ell$} [lB] at  6.355 15.334
\put{\small $z_u$} [lB] at  14.005 15.334
%
%
\linethickness= 0.500pt
\setplotsymbol ({\thinlinefont .})
{\putrule from 8.350 16.351 to 8.350 16.034
}%
%
%
\linethickness= 0.500pt
\setplotsymbol ({\thinlinefont .})
{\putrule from 14.605 16.351 to 14.605 16.034
}%
%
%
\linethickness= 0.500pt
\setplotsymbol ({\thinlinefont .})
\setdots < 0.0953cm>
{\plot  4.128 12.065 10.636 18.415 /
}%
%
%
\linethickness= 0.500pt
\setplotsymbol ({\thinlinefont .})
\setsolid
{\putrule from  4.128 12.065 to  4.128 12.067
\plot  4.128 12.067  4.130 12.069 /
\plot  4.130 12.069  4.134 12.076 /
\plot  4.134 12.076  4.138 12.086 /
\plot  4.138 12.086  4.147 12.101 /
\plot  4.147 12.101  4.157 12.120 /
\plot  4.157 12.120  4.170 12.145 /
\plot  4.170 12.145  4.187 12.177 /
\plot  4.187 12.177  4.208 12.215 /
\plot  4.208 12.215  4.231 12.262 /
\plot  4.231 12.262  4.261 12.315 /
\plot  4.261 12.315  4.293 12.376 /
\plot  4.293 12.376  4.331 12.446 /
\plot  4.331 12.446  4.373 12.522 /
\plot  4.373 12.522  4.420 12.609 /
\plot  4.420 12.609  4.470 12.702 /
\plot  4.470 12.702  4.528 12.802 /
\plot  4.528 12.802  4.587 12.912 /
\plot  4.587 12.912  4.652 13.026 /
\plot  4.652 13.026  4.722 13.151 /
\plot  4.722 13.151  4.796 13.280 /
\plot  4.796 13.280  4.875 13.415 /
\plot  4.875 13.415  4.957 13.555 /
\plot  4.957 13.555  5.044 13.701 /
\plot  5.044 13.701  5.135 13.851 /
\plot  5.135 13.851  5.228 14.008 /
\plot  5.228 14.008  5.326 14.165 /
\plot  5.326 14.165  5.425 14.326 /
\plot  5.425 14.326  5.529 14.489 /
\plot  5.529 14.489  5.635 14.654 /
\plot  5.635 14.654  5.743 14.821 /
\plot  5.743 14.821  5.855 14.988 /
\plot  5.855 14.988  5.969 15.155 /
\plot  5.969 15.155  6.083 15.323 /
\plot  6.083 15.323  6.202 15.490 /
\plot  6.202 15.490  6.322 15.655 /
\plot  6.322 15.655  6.445 15.818 /
\plot  6.445 15.818  6.568 15.981 /
\plot  6.568 15.981  6.695 16.142 /
\plot  6.695 16.142  6.824 16.298 /
\plot  6.824 16.298  6.955 16.453 /
\plot  6.955 16.453  7.087 16.603 /
\plot  7.087 16.603  7.222 16.751 /
\plot  7.222 16.751  7.360 16.895 /
\plot  7.360 16.895  7.499 17.035 /
\plot  7.499 17.035  7.641 17.170 /
\plot  7.641 17.170  7.787 17.302 /
\plot  7.787 17.302  7.935 17.429 /
\plot  7.935 17.429  8.086 17.549 /
\plot  8.086 17.549  8.240 17.666 /
\plot  8.240 17.666  8.397 17.776 /
\plot  8.397 17.776  8.558 17.879 /
\plot  8.558 17.879  8.723 17.977 /
\plot  8.723 17.977  8.890 18.068 /
\plot  8.890 18.068  9.064 18.150 /
\plot  9.064 18.150  9.239 18.227 /
\plot  9.239 18.227  9.421 18.294 /
\plot  9.421 18.294  9.605 18.354 /
\plot  9.605 18.354  9.794 18.404 /
\plot  9.794 18.404  9.986 18.445 /
\plot  9.986 18.445 10.183 18.476 /
\plot 10.183 18.476 10.384 18.495 /
\plot 10.384 18.495 10.590 18.504 /
\plot 10.590 18.504 10.797 18.502 /
\plot 10.797 18.502 11.007 18.485 /
\plot 11.007 18.485 11.218 18.457 /
\plot 11.218 18.457 11.430 18.415 /
\plot 11.430 18.415 11.625 18.364 /
\plot 11.625 18.364 11.819 18.303 /
\plot 11.819 18.303 12.012 18.231 /
\plot 12.012 18.231 12.203 18.148 /
\plot 12.203 18.148 12.391 18.057 /
\plot 12.391 18.057 12.575 17.958 /
\plot 12.575 17.958 12.757 17.852 /
\plot 12.757 17.852 12.935 17.736 /
\plot 12.935 17.736 13.109 17.615 /
\plot 13.109 17.615 13.280 17.488 /
\plot 13.280 17.488 13.447 17.355 /
\plot 13.447 17.355 13.612 17.215 /
\plot 13.612 17.215 13.771 17.071 /
\plot 13.771 17.071 13.928 16.923 /
\plot 13.928 16.923 14.080 16.768 /
\plot 14.080 16.768 14.230 16.612 /
\plot 14.230 16.612 14.376 16.449 /
\plot 14.376 16.449 14.520 16.284 /
\plot 14.520 16.284 14.660 16.114 /
\plot 14.660 16.114 14.798 15.943 /
\plot 14.798 15.943 14.933 15.767 /
\plot 14.933 15.767 15.064 15.587 /
\plot 15.064 15.587 15.196 15.405 /
\plot 15.196 15.405 15.323 15.221 /
\plot 15.323 15.221 15.447 15.035 /
\plot 15.447 15.035 15.570 14.844 /
\plot 15.570 14.844 15.693 14.652 /
\plot 15.693 14.652 15.811 14.459 /
\plot 15.811 14.459 15.928 14.262 /
\plot 15.928 14.262 16.044 14.065 /
\plot 16.044 14.065 16.159 13.866 /
\plot 16.159 13.866 16.271 13.665 /
\plot 16.271 13.665 16.381 13.464 /
\plot 16.381 13.464 16.489 13.261 /
\plot 16.489 13.261 16.595 13.060 /
\plot 16.595 13.060 16.701 12.857 /
\plot 16.701 12.857 16.802 12.656 /
\plot 16.802 12.656 16.904 12.454 /
\plot 16.904 12.454 17.003 12.253 /
\plot 17.003 12.253 17.098 12.054 /
\plot 17.098 12.054 17.194 11.860 /
\plot 17.194 11.860 17.287 11.665 /
\plot 17.287 11.665 17.376 11.472 /
\plot 17.376 11.472 17.465 11.286 /
\plot 17.465 11.286 17.549 11.102 /
\plot 17.549 11.102 17.630 10.922 /
\plot 17.630 10.922 17.710 10.746 /
\plot 17.710 10.746 17.784 10.577 /
\plot 17.784 10.577 17.858 10.414 /
\plot 17.858 10.414 17.928 10.257 /
\plot 17.928 10.257 17.994 10.107 /
\plot 17.994 10.107 18.055  9.963 /
\plot 18.055  9.963 18.114  9.828 /
\plot 18.114  9.828 18.169  9.701 /
\plot 18.169  9.701 18.220  9.580 /
\plot 18.220  9.580 18.267  9.468 /
\plot 18.267  9.468 18.311  9.364 /
\plot 18.311  9.364 18.351  9.271 /
\plot 18.351  9.271 18.387  9.184 /
\plot 18.387  9.184 18.419  9.106 /
\plot 18.419  9.106 18.449  9.038 /
\plot 18.449  9.038 18.474  8.977 /
\plot 18.474  8.977 18.495  8.924 /
\plot 18.495  8.924 18.514  8.877 /
\plot 18.514  8.877 18.529  8.841 /
\plot 18.529  8.841 18.542  8.810 /
\plot 18.542  8.810 18.553  8.784 /
\plot 18.553  8.784 18.561  8.765 /
\plot 18.561  8.765 18.565  8.750 /
\plot 18.565  8.750 18.570  8.742 /
\plot 18.570  8.742 18.572  8.735 /
\plot 18.572  8.735 18.574  8.733 /
\putrule from 18.574  8.733 to 18.574  8.731
}%
\linethickness=0pt
\putrectangle corners at  3.150 23.520 and 23.203  7.595
\endpicture}
\end{center}
\nobreak
\vspace*{-0.5cm}
\begin{center}
  {\small\bf Figure~1}
\end{center}

For $t\in \N$ let $z(t)=1+\sqrt{2}\|x[t^*]\|$. 
When $t=0$ we have $t^*=0$ as well and 
$$
  z(0)\leq 
  1+\sqrt{2}\|x[0]\|
  =\bb< \left(\frac1{\ba}\right)^{\frac1{\a-1}}
      \left(\frac1{\a}\right)^{\frac1{\a-1}}=z_*.
$$
This implies that $z(0)< z_\ell$. 
Assume that there exists $t<T$ such that $z(t)\geq z_u$ and 
let $r$ be the first such $t$. Then $r=r^*\geq 1$ and, for 
all $t<r$ 
$$
  1+\sqrt{2}\|x[t]\| \leq z(r-1)\leq z_\ell.
$$
Let $z_0$ be the intersection of the $z$ axis with the line 
segment joining $(0,-\bb)$ and $(z_*,F(z_*))$ (see Figure~1).  
The line where this segment lies has equation 
$$
   y+\bb=z\frac{z_*-\ba z_*^\a}{z_*}
$$
from which it follows that 
$$
  z_0=\frac{\bb}{1-\ba z_*^{\a-1}}
     =(1+\sqrt{2}\|x(0)\|)\frac{\a}{\a-1}.
$$
It follows that, for all $t<r$, 
$$
  \|x[t]\|\leq 
  \frac{1}{\sqrt{2}}\left(z_\ell-1\right)
  \leq \frac{1}{\sqrt{2}}\left(z_0-1\right)=B_0.
$$
In particular,
$$
  \|x[r-1]\|\leq \frac{1}{\sqrt{2}}(z_\ell-1).
$$
For $r$ instead, we have
$$
  \|x[r]\|\geq \frac{1}{\sqrt{2}}(z_u-1).
$$
This implies 
\begin{equation}\label{eq:sup}
  \|x[r]-x[r-1]\| \geq 
   \|x[r]\|-\|x[r-1]\| \geq 
   \frac{1}{\sqrt{2}}(z_u-z_\ell)\geq 
   \frac{1}{\sqrt{2}}(z_*-z_\ell).
\end{equation}
From the intermediate value theorem, there is $\xi\in[z_\ell,z_*]$ 
such that $F(z_*)=F'(\xi)(z_*-z_\ell)$. But $F'(\xi)\geq 0$ and 
$F'(\xi)=1-a\a\xi^{\a-1}\leq 1$. 
Therefore, 
$$
   z_*-z_\ell\geq F(z_*)
$$
and it follows from (\ref{eq:sup}) that 
\begin{equation}\label{eq:sup2}
  \|x[r]-x[r-1]\| \geq \frac{1}{\sqrt{2}}F(z_*).
\end{equation}
But
$$
\|x[r]-x[r-1]\| = h\|v[r-1]\| \leq h \|v[0]\|
$$
the last since $\|v\|$ is decreasing for $t<T$.  
Putting this inequality together with (\ref{eq:sup2}) shows that  
$$
  F(z_*)\leq \sqrt{2}h\|v[0]\|
$$
or equivalently,
$$
 \left(\frac1{\ba}\right)^{\frac1{\a-1}}
         \left(\frac1{\a}\right)^{\frac1{\a-1}}
         \frac{\a-1}{\a}-\bb
  \leq \sqrt{2}h\|v[0]\|
$$
which contradicts our hypothesis. This shows that, for all 
$t<T$, $z(t)< z_\ell$ and hence, for all $t<\hat t$, 
$$
  (1+\sqrt{2}\|x[t]\|)^\a \leq z_0^\a
  =\left((1+\sqrt{2}\|x[0]\|)\frac{\a}{\a-1}\right)^\a
  \leq 2^{-\a}\frac{kK}{2\HH_0}
$$ 
the last by the definition of $\HH_0$. 
We now proceed as in case~(i). 
\eproof

\subsection{Proof of Theorem~\ref{th:disc_absolute}}
\label{sec:disc_abs}

\begin{proposition}
For $\e>0$, let $p(\e)=\Prob\{\|H\|\leq \e\}$. Then, in the
uniform case, we have the bound
$$
   p(\e)\geq\left(\frac{\e}{r}\right)^{3k},
$$
($p(\e)=1$ if ${r}\leq \e$) 
while in the Gaussian case $\|H/\sigma\|^2$ has a 
Chi-square distribution
with $3k-3$ degrees of freedom, and in consequence
$$
p(\e)=\int_0^{\sqrt{\e/(2\sigma)}}\frac{t^{\frac{3k-5}{2}}}
 {\Gamma\left(\frac{3k-3}{2}\right)}e^{-t}dt.
$$
\end{proposition}

\proof 
In the uniform case, since $\|H\|\leq\|\bH\|$, we have
\begin{equation*}
  \Prob(\|H\|\leq\e)\geq
  \Prob(\|\bH\|\leq\e)
  =\left(\frac{\e}{r}\right)^{3k}.
\end{equation*}
In the Gaussian case, the decomposition
$\bH=\bH_{\Delta}+\bH_{\perp}$ takes the form
$\bH_{\Delta}=(\bm,\ldots,\bm)$, with 
$\bm=\frac 1k\sum_{j=1}^k\bH_j$,
and $\bH_{\perp}=\bH-(\bm,\ldots,\bm)$. Consequently
\begin{align*}
 \|H\|^2=\|\bH_{\perp}\|^2=\sum_{j=1}^k(\bH_j-\bm)^2
 =\sum_{\ell=1}^3\sum_{j=1}^k(e^{(\ell)}_j-\bm^{(\ell)})^2,
\end{align*}
where $\bm=\left(\bm^{(1)},\bm^{(2)},\bm^{(3)}\right)$. 
A standard result in
statistics (see \cite{Arnold90}, page 219) states that
$\frac 1{\sigma^2}\sum_{j=1}^k(e^{(\ell)}_j-\bm^{(\ell)})^2$ has a 
Chi-square distribution with $k-1$ degrees of freedom. 
Therefore, by independence, $\|H/\sigma\|^2$ has a
Chi-square distribution with $3k-3$ degrees of freedom. The 
expression for $p(\e)$ follows from the form of the density 
of the Chi-square random variable. 
\eproof

We can now give the proof of Theorem~\ref{th:disc_absolute}.

Let $T>0$ and assume that the hypothesis of Theorem~\ref{th:2}
holds for $T$. Then,
$$
 \|v[t]\|\leq \|v[0]\|
 \left(1-h\frac{kK}{2(1+\sqrt{2}B_0)^\a}\right)^t
$$
where $B_0$ is as in Theorem~\ref{th:2}. Therefore
\begin{eqnarray*}
\|v[T]\|\leq \nu &\Longleftarrow &
  \|v[0]\|\left(1-h\frac{kK}{2(1+\sqrt{2}B_0)^\a}\right)^T \leq \nu\\
&\iff & T\geq \left(\ln
\left(1-h\frac{kK}{2(1+\sqrt{2}B_0)^\a}\right)\right)^{-1}
 \ln\left(\frac{\nu}{\|v[0]\|}\right)\\
&\Longleftarrow & T\geq
\frac{2(1+\sqrt{2}B_0)^\a}{hkK}\ln\left(\frac{\|v[0]\|}{\nu}\right)
  =T_0.
\end{eqnarray*}
At each iteration $t$, Theorem~\ref{th:2} requires that
$\|H(t)\|\leq \HH_0\|v[t]\|$ for a quantity $\HH_0$ which depends on
the initial conditions and on the case ((i), (ii), or (iii))
at hand. If nearly-alignment has not occurred, then
$\|v[t]\|\geq \nu$ and therefore,
$$
 \Prob\{\|H(t)\|\leq \HH_0\|v[t]\|\} \geq
 \Prob\{\|H(t)\|\leq \HH_0\nu\} \geq p(\HH_0\nu).
$$
It follows that
$$
   \Prob\{\|H(t)\|\leq \HH_0\|v[t]\|
    \mbox{ for $t=0,\ldots,T_0-1$}\} \geq p(\HH_0\nu)^{T_0}
$$
and therefore the claimed bounds in the uniform and Gaussian cases. 
\eproof

\section{Continuous time}\label{sec:cont_time}

The goal of this section is to show a continuous time 
version of Theorem~\ref{th:disc_absolute}. In contrast with the 
discrete time setting, though, the description of the noise 
$\bH$ is not straightforward. There is no obvious continuous time version of 
a sequence of independent and identically distributed 
random variables. We thus begin by 
discussing our model of noise. 

\subsection{Continuous time stochastic processes}\label{sec:CSP}

A {\em continous time stochastic process} $\{X(t)\mid  t\geq 0\}$, or
for short a {\em stochastic process}, is a family of
random variables $X(t)$ defined in a common probability space
$(\Omega,\FF,\Prob)$. More precisely, the stochastic process depends 
on two arguments, $t$ and $\omega$, and
when necessary we denote it by $X(t,\omega)$. In order to
describe our assumptions, we make a brief comparison with the
discrete time situation.

In the discrete time situation it is natural to assume that each
perturbation $X(t)$ is a centered random variable 
(i.e., $\bfE(X(t))=0$), that for all values of $t$
the random variables $X(t)$ are mutually
independent, and that they have a common distribution. 
Making some additional assumption on the distribution of each
random variable, typically assuming normality with a given
standard deviation $\sigma$, completes the specification of the
noise probabilistic structure. Let us call this structure a
{\em Gaussian white noise sequence}.

Unfortunately, in continuous time we do not have a reasonable
tractable model that shares the properties of the Gaussian white
noise sequence. To understand why this is so we need to consider
the properties of the {\em trajectories} of the processes, i.e., 
the curves
obtained when $\omega\in\Omega$ is fixed and $t$ ranges in the
interval $[0,\infty)$. Assuming independence
of the random variables of the stochastic process for
close values of $t$ makes the trajectories of the process to have
an extremely irregular behaviour. Consequently,
if we want the trajectories of the process to be continuous, or
differentiable, we can not assume the independence of the random
variables for pairs of close values of $t$. 

To analyze this difficulty consider 
first, in the discrete time case, the accumulated perturbation
produced by a Gaussian white noise sequence
$\{X(t)\mid t\in\N\}$. 
That is, consider the
sequence of sums
\begin{equation*}
  S(0)=0,\qquad S(t)=X(1)+\cdots+X(t),\quad t=1,2,\ldots.
\end{equation*}
This random sequence associated to the Gaussian white noise
sequence, called {\em Gaussian random walk}, does have a natural
counterpart in the continuous time case, known as {\em Wiener
process} or {\em Brownian motion}. A Wiener process 
$W=\{W(t)\mid t\geq 0\}$ is a continuous time stochastic 
process satisfying 
the following properties, that are natural extensions
to continuous time of the properties of the sequence 
$\{S(t)\mid t\in\N\}$: 
\begin{description}
\item{(a)}
$W(0)=0$.
\item{(b)}
The increments of the process are
independent random variables. That is, for all
$0\leq t_0<t_1<\dots<t_n$ the random variables
\begin{equation*}
    W(t_1)-W(t_0),\dots,W(t_n)-W(t_{n-1}),
\end{equation*}
are independent.
\item{(c)} The increments are homogeneous and
Gaussian. That is, for all $t,h\geq 0$,  the random
variable
\begin{equation*}
   W(t+h)-W(t)
\end{equation*}
has a centered Gaussian distribution with variance $h$.
\item{(d)}
Finally, the trajectories of $W$ are continuous. That is, the
curves $\{W(\cdot,\omega)\}$ obtained when $\omega\in\Omega$ is
fixed and $t\geq 0$ are continuous functions of $t$ for
almost all $\omega\in\Omega$.
\end{description}
In the discrete time case the
Gaussian white noise sequence can be recovered from the sums
$S(t)$ by taking differences $X(t)=S(t)-S(t-1)$. In the
continuous time case we would like to have a formula like
$X(t)=\dot{W}(t)$, but it is not immediate to give a sense to this
last time derivative, as it is known that the trajectories of the
Wiener process are nowhere differentiable~\cite{Billingsley95}. 
A possible way out is to first take a differentiable approximation
of the Wiener process and then take the time derivative
of this approximation as a model of noise in the continuous time
case. Another, alternative, way out 
is described in Remark~\ref{remark:ito} below.

The approximation is obtained by convolution with a smooth
kernel. Let $\psi\colon\R\rightarrow\R_{+}$ be a 
$\CC^{1}$ function, 
with compact support, say  
$\textsf{supp}(\psi)\subset [-{1/2},{1/2}]$
and such that $\int_{\R}\psi(x)dx=1$.
For $\delta>0$ consider
$$
  \psid(x)=\frac1{\delta}\psi\Big(\frac{x}{\delta}\Big),
$$
that has $\textsf{supp}(\psid)\subset [-{\delta/2},{\delta/2}]$.

The approximation
$W^{\delta}=\{W^{\delta}(t)\mid t\geq 0\}$ of the Wiener
process is obtained by convolution with $\psid$ in the
following way:
\begin{equation}
  W^{\delta}(t)=\big(\psid \ast W\big)(t)
  =\int_{\R}\psid(t-s)W(s)ds
  =\int_{\R}\psid(-w)W(t+w)dw, \label{reg}
\end{equation}
where, if $s<0$,  we replace $W(s)$ in the integrand by 
$\hat{W}(-s)$, where
$\hat{W}$ is another Wiener process independent of $W$.

Observe that the process $W^\delta$ inherits the regularity
properties of $\psi$. In particular, it has $\CC^1$ 
trajectories. 
We now define a {\em noise process} $X^{\delta}$ as the time derivative 
of $W^\delta$,
\begin{align}
  X^{\delta}(t)&=\frac{d}{dt} W^{\delta}(t)
    =\int_{\R}{\frac{\partial }{\partial t}}
    \big(\psid(t-s)\big)W(s)ds\label{eq:1}\\
  &=\frac1{\delta}\int_{\R}\dot{\psi}(-w)
    \left(W(t+\delta w)-W(t-{\delta/2})\right) dw  \label{eq:2}\\
  &=\int_{\R}\psid(t-s)dW(s).\label{eq:3}
\end{align}
Here the second equality in \eqref{eq:1} is obtained 
by differentiation under the integral sign, 
\eqref{eq:2} after a change of variables $s-t=\delta w$ and using that
$\int_{\R}\dot{\psi}(-w)dw=0$, and   
\eqref{eq:3} by integration by parts departing from 
\eqref{eq:1} (note, this involves a 
stochastic integral \cite{Oksendal03}). 

Using~\eqref{eq:3} and It\^o's isometry for the stochastic 
integral~\cite[Theorem 4.2]{Oksendal03} 
we obtain (we write $X$ instead of $X^\delta$), for $h\geq 0$, 
\begin{align}\notag
  \bfE\left(X(t+h)X(t)\right)
  &=\bfE\left(\int_{\R}\psid(t+h-s)dW(s)
          \int_{\R}\psid(t-s)dW(s)\right)\\
  &=\int_{\R}\psid(t+h-s)\psid(t-s)ds.\label{eq:isometry}
\end{align}
Taking $h=0$, we obtain
\begin{equation}\label{eq:consequences}
  \bfE\left(X(t)^2\right)
  =\int_{\R}\psid(t-s)\psid(t-s)ds=\|\psid\|^2
  =\frac1{\delta}\|\psi\|^2.
\end{equation}
Taking $h\geq\delta$ in \eqref{eq:isometry} and using that 
$\textsf{supp}(\psid)\subset [-{\delta/2},{\delta/2}]$ we obtain 
\begin{equation}\label{eq:consequences2}
   \bfE\left(X(t+h)X(t)\right) =0.
\end{equation}
Furthermore, from the different expressions
in (\ref{eq:1}--\ref{eq:3}) we obtain the following
properties of the process $\{X(t)\mid t\geq 0\}$:
\begin{description}
\item{(a)}
It is a {\em centered Gaussian} process with variance
$\|\psid\|^2=(1/\delta)\|\psi\|^2$. 
That is, $X(t)\sim N(0,(1/\delta)\|\psi\|^2)$ for all $t\in[0,T]$.
\item{(b)}
It is a {\em stationary} process. That is, for all times
$t_1,\dots,t_n$, intervals $I_1,\dots,I_n$, and 
time increment $h$, we have
\begin{equation*}
  \Prob(X(t_1)\in I_1,\dots,X(t_n)\in I_n)=\Prob(X(t_1+h)\in
  I_1,\dots,X(t_n+h)\in I_n).
\end{equation*}
\item{(c)}
It is {\em $\delta$-dependent}. That is, the two sets of random
variables
\begin{equation*}
  \{X(s)\mid s\leq t\}\quad\text{and}\quad
  \{X(s)\mid s\geq t+\delta\}
\end{equation*}
are independent for each $t$. 
\end{description}
The first property can be obtained from \eqref{eq:1} 
since the integral there is the limit of a linear combination 
of Gaussian random variables. Indeed, such a linear combination 
remains Gaussian and the limit of the resulting random variables 
preserves Gaussianity as well. This variable is centered 
since all the involved variables are centered, and the 
limit defining the
integral preserves the expectation. The value for the variance 
follows from \eqref{eq:consequences}. 

The stationarity property of $X(t)$ is
inherited from the stationarity of the increments of the Wiener
process. In order to see it we use the representation
\eqref{eq:2}. 
The probability distribution of the stochastic process
$\{W(t+w)-W(t-\delta/2)\mid -\delta/2\leq w\leq \delta/2\}$
does not depend on the value of $t$, or more precisely, 
the probability distribution of 
$\{W(t_i+w)-W(t_i-\delta/2)\mid -\delta/2\leq w\leq \delta/2,
i=1,\dots,n\}$
coincides with the probability distribution of
$\{W(t_i+h+w)-W(t_i+h-\delta/2)\mid -\delta/2\leq w\leq \delta/2,
i=1,\dots,n\}$.
This makes the process probabilities invariant under
a shift of $h$, i.e. the process satisfy the definition 
of stationarity in~(b) above. 

Finally, the $\delta$-dependency is a consequence of 
\eqref{eq:consequences2}. This equality give us non-correlation,
when the lag $h\geq\delta$. The independence follows since 
a Gaussian vector without correlation has independent components.
It should be noticed that this property is not essential to our
developement below; it simply mimics the discrete time 
independence.
Furthermore, it is possible to derive similar results in this
discrete time case for a weakly dependent noise.

\begin{remark}
The contents of this section is not new. It is exposed in 
certain detail for ease of the reader. Regarding the equalities in 
(\ref{eq:1}--\ref{eq:3}), it can be seen that any 
centered Gaussian process
admits such a representation with an adequate kernel.
General results on Gaussian processes and their diverse
applications can be found for instance 
in~\cite{Adler90,aw,CL:2004}.
\end{remark}

\subsection{Statement of the main result}\label{sec:stat_cont}

We can now describe the noise $\bH$ in \eqref{eq:dif_cont} 
and state a continuous version of Theorem~\ref{th:disc_absolute}. 

We assume that $\bH_i(t)$ is a three dimensional Gaussian
centered, stationary stochastic process, that satisfies a
$\delta$-dependence condition for some $\delta>0$, has $\CC^1$ 
trajectories, and independent coordinates. 
More precisely, we assume that 
$\bH_i(t)=\left(e_i^{(1)}(t),e_i^{(2)}(t), e_i^{(3)}(t)\right)$, 
where each coordinate is given by 
\begin{equation}\label{eq:representation}
   e_i^{(\ell)}(t)=\sigma\sqrt{\delta}
       \int_{\R}\psid(t-s)dW_i^{(\ell)}(s)
\end{equation}
where $\psid$ is a kernel as in \S\ref{sec:CSP}, $\sigma>0$, 
and $\{W_i^{(\ell)}(t)\mid t\geq 0\}$ is a set of $3k$ independent
Wiener processes. That is, each coordinate of $\bH_i$ is of the 
form $\sigma\sqrt{\delta}X^{\delta}$ with $X^{\delta}$ as 
in~\S\ref{sec:CSP}. Note that the variance 
$\Var\left(e_i^{(\ell)}(t)\right)=\sigma^2$ for all $t\geq 0$. 

\begin{theorem}\label{th:main_cont}
Consider the system \eqref{eq:dif_cont} with adjacency 
matrix given by~\eqref{eq:adj} and noise given 
by~\eqref{eq:representation}.
Let $x_0,v_0\in\E$. Then, there exists 
a unique solution 
$(x(t),v(t))$ of \eqref{eq:dif_cont}, defined for all $t\in\R$,  
with initial conditions $x(0)=x_0$ and $v(0)=v_0$. 
Assume that one of the three following hypothesis 
holds:
\begin{description}
\item[(i)]
$\a<1$, 
\item[(ii)]
$\a=1$, and $\|v(0)\|<\frac{kK}{2\sqrt{2}}$, 
\item[(iii)]
$\a>1$,  and 
$$ 
  \left(\frac1{\a\ba}\right)^{\frac1{\a-1}}
  \frac{\a-1}{\a}>\bb. 
$$
\end{description} 
Then $\nu$-nearly-alignment for some $\nu<\|v(0)\|$ occurs before time 
$$
  T_0:=\frac{U_0^\a}{kK}
  \ln\left(\frac{\|v(0)\|}{\nu}\right)
$$
with probability at least
\begin{equation}\label{eq:probability}
 \left\{
  2\bPhi\left(\frac{\nu\HH_0}{\sigma\sqrt{3k}}\right)
  -\frac{2T_0\sigma\|\dot{\psi}\|}{\delta\sqrt{2\pi}} 
  \varphi\left(\frac{\nu\HH_0}{\sigma\sqrt{3k}}\right)-1
 \right\}^{3k}
\end{equation}
where $\bPhi(y)=\frac{1}{\sqrt{2\pi}}\int_{-\infty}^y 
e^{-\frac{u^2}{2}} du$ is the standard normal distribution
and $ \varphi=\bPhi'$ its density.
\end{theorem}

\begin{remark}\label{rem:cont}
Using the identity
$$
 1-\bPhi(x)=\frac{\varphi(x)}{x}
  \left(1+\Oh\left(\frac 1{x^2}\right)\right)
$$
and performing some elementary computations 
we obtain that the bound in \eqref{eq:probability} 
is equivalent, for small $\sigma$, to
$$
 1-6k\sigma\left(\frac{\sqrt{3k}}{\nu\HH_0}+
 \frac{T_0\|\dot{\psi}\|}{\delta\sqrt{2\pi}}\right)
 \varphi\left(\frac{\nu\HH_0}{\sigma\sqrt{3k}}\right).
$$
\end{remark}

\begin{remark}\label{remark:ito}
An alternative way to model the noise in our
system relies on the similarity of the Gaussian random walk
and the Wiener process. Integrating the second equation 
in~\eqref{eq:dif_cont} we obtain
$$
  \bv(t)=\bv(0)-\int_0^tL_x\bv(s)ds +\int_0^t\bH(s)ds.
$$
The last term in the right-hand side is the accumulated noise 
for which we noted in \S\ref{sec:CSP} that the natural   
continous time version is the Wiener process. Multiplying 
the latter by $\sigma>0$ (as we did at the beginning of 
\S\ref{sec:stat_cont} to obtain $\Var(\bH)=\sigma^2\Id_{3k}$) 
we obtain 
$$
  \bv(t)=\bv(0)-\int_0^tL_x\bv(s)ds +\sigma W(t)
$$
and integral equation often written in its ``differential form'' 
$$
   d\bv(t) = - L_x\bv(t)dt +\sigma d W(t).
$$
Hence, an alternative to~\eqref{eq:dif_cont} is the system of 
stochastic differential equations
\begin{align}\label{eq:SDE}
d\bx\;&=\; \bv dt\tag{SDE}\\ 
d\bv\;&=\; -L_x \bv dt+ \sigma dW.\notag
\end{align}
The construction of a solution for this system relies on the 
{\em stochastic calculus} developed by It\^o~\cite{Ito06}. 

We note that, while it is possible  to prove that $W^{\delta}\to W$
when $\delta\to0$, it is not generally true (cf.~\cite{KP:91}) 
that the solution of a
system of stochastic differential equations driven by a smoothed noise
$W^{\delta}$ converges towards the solution of the corresponding system
driven by the original noise $W$. 
Investigating whether this is the case for~\eqref{eq:dif_cont} 
and~\eqref{eq:SDE} would take us out of the scope of the present 
work.
\end{remark}

To prove Theorem~\ref{th:main_cont} we follow the steps in the 
proof of Theorem~\ref{th:disc_absolute}. 

\subsection{Bounded noise}\label{sec:cont_rob}

For $x\in X$ we denote $\Gamma(x)=\|x\|^2$ and for 
$v\in V$ we denote $\Lambda(v)=\|v\|^2$.  

In this section we fix $T\in(0,\infty)\cup\{+\infty\}$ 
and a solution $(x,v)$ of (\ref{eq:dif_cont}) 
(which we assume exists and is, almost surely, 
differentiable in $[0,T)$).  
The meaning of expressions like $\phi_t$, 
$L_t$, $\Lambda(t)$, or $\Gamma(t)$ 
is as described in \S\ref{sec:disc_rob}.  

Denote $\displaystyle\Phi_t=\min_{\tau\in[0,t]}\phi_{\tau}$. 

\begin{proposition}\label{prop:decay}
Assume that, for all $0\leq t<T$, 
$\|H_t\|\leq \|v(t)\|\HH_0$. Then, for all $0\leq t<T$,  
$$
  \|v(t)\|\leq \|v(0)\|e^{-t(\Phi_t-\HH_0)}. 
$$
\end{proposition}

\proof
Let $\tau\in[0,t]$. Then 
\begin{eqnarray*}
 \Lambda'(\tau)
&=&\frac{d}{d\tau}\langle v(\tau),v(\tau)\rangle\\ 
&=&2\langle v'(\tau),v(\tau)\rangle\\ 
&=&-2\langle L_{\tau}v(\tau)+H_{\tau},v(\tau)\rangle\\ 
&=&-2\langle L_{\tau}v(\tau),v(\tau)\rangle
   -2\langle H_{\tau},v(\tau)\rangle\\ 
&\leq &-2\phi_{x(\tau)}\Lambda(\tau)+
    2\|H_{\tau}\|\|v(\tau)\|\\
&\leq &-2\Lambda(\tau)(\phi_{x(\tau)}-\HH_0).
\end{eqnarray*}
Here we used that $L_{\tau}$ is symmetric positive semidefinite 
on $V$. Using this inequality, 
$$
  \ln(\Lambda(\tau))\biggl|_0^t 
  =\int_0^t\frac{\Lambda'(\tau)}{\Lambda(\tau)}d\tau 
  \leq \int_0^t -2(\phi_{\tau}-\HH_0) d\tau 
  \leq -2t(\Phi_t-\HH_0)
$$
i.e.,
$$
  \ln(\Lambda(t))-\ln(\Lambda_0)\leq -2t(\Phi_t-\HH_0)
$$
from which the statement follows.
\eproof

\begin{proposition}\label{prop:T}
Assume that $\Phi_t>\HH_0$ for all $0\leq t<T$. Then, 
for all $0\leq t<T$,
$$
  \|x(t)\|\leq \|x(0)\|+
  \frac{\|v(0)\|}{\Phi_t-\HH_0}.
$$
\end{proposition}

\proof
For $\tau\leq t$ we have 
$|\Gamma'(\tau)|=|2\langle v(\tau),x(\tau)\rangle|\leq 2\|v(\tau)\| 
\|x(\tau)\|$. But $\|x(\tau)\|=\Gamma(\tau)^{1/2}$ and  
$\|v(\tau)\|^2=\Lambda(\tau)\leq 
\Lambda_0e^{-2\tau(\Phi_\tau-\HH_0)}$, by 
Proposition~\ref{prop:decay}. Therefore, 
\begin{equation}\label{eq:gamma'}
  \Gamma'(\tau)\leq |\Gamma'(\tau)|\leq 
  2\left(\Lambda_0e^{-2\tau(\Phi_\tau-\HH_0)}\right)^{1/2}
  \Gamma(\tau)^{1/2}
\end{equation}
and, using that $\tau\mapsto \Phi_\tau$ is non-increasing and that 
$\Phi_\tau-\HH_0>0$ for all $\tau\leq t$, 
\begin{eqnarray*}
  \int_0^t\frac{\Gamma'(\tau)}{\Gamma(\tau)^{1/2}}d\tau
  &\leq& 2\int_0^t 
  \left(\Lambda_0e^{-2\tau(\Phi_\tau-\HH_0)}\right)^{1/2} d\tau\\
  &\leq& 2\int_0^t \Lambda_0^{1/2}e^{-\tau(\Phi_t-\HH_0)} d\tau\\
  &=& 2\Lambda_0^{1/2}\left(-\frac{1}{\Phi_t-\HH_0}\right)
      e^{-\tau(\Phi_t-\HH_0)}\biggl|_0^t 
  \;\leq\; \frac{2\Lambda_0^{1/2}}{\Phi_t-\HH_0}
\end{eqnarray*}
the last inequality because $\Phi_t>\HH_0$. This implies
$$
    \Gamma(\tau)^{1/2}\biggl|_0^t
    = \frac12\int_0^t\frac{\Gamma'(\tau)}{\Gamma(\tau)^{1/2}}d\tau 
    \leq \frac{\Lambda_0^{1/2}}{\Phi_t-\HH_0}
$$
from which it follows that
\begin{equation}\tag*{\qed}
  \Gamma(t)^{1/2}\leq \Gamma_0^{1/2}+
  \frac{\Lambda_0^{1/2}}{\Phi_t-\HH_0}.
\end{equation}

The main result in this section is the following. 

\begin{theorem}\label{th:1}
Assume that, for 
all $0\leq t<T$, $\|H(t)\|\leq \|v(t)\|\HH_0$.   
Assume also that one of the three following hypothesis 
hold: 
\begin{description}
\item[(i)]
$\a<1$, 
\item[(ii)]
$\a=1$, and $\|v(0)\|<\frac{kK}{2\sqrt{2}}$,  
\item[(iii)]
$\a>1$, and 
$$ 
  \left(\a\frac1{\ba}\right)^{\frac1{\a-1}}
  \frac{\a-1}{\a}>\bb. 
$$
\end{description}
Then, for all $0\leq t <T$, $\|x(t)\|\leq B_0$ and 
$$
  \|v(t)\|\leq \|v(0)\| e^{-\frac{kK}{U_0^\a}t}.
$$
In particular, when $T=\infty$, $\|v(t)\|\to0$ for $t\to\infty$ 
and there exists $\hat{x}\in X$ such that 
$x(t)\to\hat{x}$ when $t\to\infty$. 
\end{theorem}

\proof
Let 
$$
   \Upsilon=\left\{t\in [0,T)\mid (1+\sqrt{2}\|x(t)\|)^\alpha 
       \leq \frac{kK}{2\HH_0}\right\}. 
$$
Note that in all three cases ((i), (ii), and (iii)) the definition 
of $\HH_0$ implies that $0\in\Upsilon$ and hence, that 
$\Upsilon\neq\emptyset$. Assume that $\Upsilon\neq [0,T)$ and 
let $\hat t=\inf \{[0,T)\setminus\Upsilon\}$. Clearly, 
$1+\sqrt{2}\|x(\hat{t})\| = \frac{kK}{2\HH_0}$.  

By Proposition~\ref{prop:phi} and Lemma~\ref{lem:gamma}, 
for all $x\in X$, 
$$
  \phi_x\geq\frac{kK}{(1+\max_{i\neq j}\|x_i-x_j\|)^\a}
  \geq\frac{kK}{(1+\sqrt{2}\|x\|)^\a}.
$$ 
Let $t<\hat t$ and $t^*\in[0,t]$ be the point maximizing 
$\|x\|$ in $[0,t]$. Then  
$$
  \Phi_{t}=\min_{\tau\in[0,t]}\phi_{\tau} 
         \geq \min_{\tau\in[0,t]}
        \frac{kK}{(1+\sqrt{2}\|x(\tau)\|)^\a}  
        \geq \frac{kK}{(1+\sqrt{2}\|x(t^*)\|)^\a}.
$$
Moreover, since $t^*\leq t<\hat t$, $t^*\in\Upsilon$  
and we have 
\begin{equation}\label{eq:error}
  \Phi_{t}-\HH_0
  \geq \frac{kK}{(1+\sqrt{2}\|x(t^*)\|)^\a}-\HH_0
  \geq \frac{kK}{2(1+\sqrt{2}\|x(t^*)\|)^\a}>0.  
\end{equation}
Hence, we may apply Proposition~\ref{prop:T} to obtain 
\begin{eqnarray}\label{eq:bound}
  \|x(t)\|&\leq& \|x(0)\|+\|v(0)\|
  \frac{1}{\Phi_t-\HH_0}\notag\\
  &\leq& \|x(0)\|+\frac{2\|v(0)\|(1+\sqrt{2}\|x(t^*)\|)^{\a}}
  {kK}.
\end{eqnarray}
Since $t^*$ maximizes $\Gamma$ in $[0,t]$ it also does so 
in $[0,t^*]$. Thus, for $t=t^*$, (\ref{eq:bound}) takes the 
form
\begin{equation}\label{eq:t*}
  \left(1+\sqrt{2}\|x(t^*)\|\right)-2\sqrt{2}\|v(0)\|
  \frac{(1+\sqrt{2}\|x(t^*)\|)^{\a}}{kK}
  -\left(1+\sqrt{2}\|x(0)\|\right)\leq 0.
\end{equation}
Let $z=1+\sqrt{2}\|x(t^*)\|$. 
Then (\ref{eq:t*}) can be rewritten as $F(z)\leq 0$ 
with $F(z)=z-\ba z^{\a}-\bb$. 
One can now finish the proof by dividing in cases 
as in Theorem~\ref{th:1} and following the steps in its proof.
\eproof

\subsection{Proof of Theorem~\ref{th:main_cont}}

We begin with a result on the behaviour of the maximum
of the processes $e_i^{(\ell)}$ described 
in~\S\ref{sec:stat_cont}. 

\begin{proposition}\label{prop:inequality} 
Fix $T>0$.  Denote $e(t)=e_i^{(\ell)}(t)$ and
\begin{equation*}
   p(x)=\Prob\left(\max_{0\leq t\leq T}|e(t)|<x\right).
\end{equation*}
Then
$$
  1-p(x)=\Prob\left(\max_{0\leq t\leq T}|e(t)|\geq x\right)
  \leq
  2\left[\frac{T\sigma\|\dot{\psi}\|}{\delta\sqrt{2\pi}} 
  \varphi\left(\frac{x}{\sigma}\right)+
  1-\bPhi\left(\frac{x}{\sigma}\right)\right].
$$
\end{proposition}

\proof The proof is an application of Davies's 
inequality~\cite{davies:1977}
(see also Chapter~4 in~\cite{aw}):
$$
  \Prob\left(\max_{0\leq t\leq T}e(t)\geq x\right)\leq 
  \frac1{\sqrt{2\pi}}\varphi\Big(\frac{x}{\sigma}\Big)\int_0^T
  \sqrt{r_{11}(t,t)}dt+1-\bPhi\Big(\frac{x}{\sigma}\Big).
$$
where 
$
  r_{11}(t,t)=\frac{\partial^2}{\partial t^2}
   \bfE\left(e(t)^2\right)
$
and therefore, using~\eqref{eq:isometry}, 
$$
  r_{11}(t,t)=\frac{\sigma^2}{\delta^2}\|\psip\|^2. 
$$
The conclusion now follows from the trivial bound
\begin{equation}\tag*{\qed}
 \Prob\left(\max_{0\leq t\leq T}|e(t)|\geq x\right)\leq 
 2\Prob\left(\max_{0\leq t\leq T}e(t)\geq x\right).
\end{equation}
\medskip

We can now give the proof of Theorem~\ref{th:main_cont}. 
The existence of a unique solution follows, for each 
$\omega\in\Omega$, from~\cite[Chapter~8]{HiSm74}. 

Using that $\|\bH(t)\|\leq \sqrt{3k}\|\bH(t)\|_\infty$,  
and Proposition~\ref{prop:inequality} 
for one coordinate we obtain 
\begin{align*}
 \Prob\left(\max_{0\leq t\leq T}\|H\left(t\right)\|{<}\e\right)
 &\geq\Prob\left(\max_{0\leq t\leq T}
  \|\bH\left(t\right)\|{<}\e\right)\\
&\geq \Prob\left(\max_{0\leq t\leq T} 
 \max_{1\leq j\leq k}\max_{1\leq \ell\leq 3}
 |e^{(\ell)}_i\left(t\right)|{<}
  \frac{\e}{\sqrt{3k}}\right)\\
&=\Prob\left(\max_{0\leq t\leq T}|e\left(t\right)|{<}
\frac{\e}{\sqrt{3k}}\right)^{3k}\\
&\geq 
\left\{1-2\left[\frac{T\sigma\|\dot{\psi}\|}{\delta\sqrt{2\pi}} 
\varphi\left(\frac{\e}{\sigma\sqrt{3k}}\right)+
1-\bPhi\left(\frac{\e}{\sigma\sqrt{3k}}\right)\right]\right\}^{3k}\\
&=
\left\{
2\bPhi\left(\frac{\e}{\sigma\sqrt{3k}}\right)
-\frac{2T\sigma\|\dot{\psi}\|}{\delta\sqrt{2\pi}} 
\varphi\left(\frac{\e}{\sigma\sqrt{3k}}\right)-1
\right\}^{3k}.
\end{align*}
Similarly as in the proof of Theorem~\ref{th:disc_absolute},
but taking into account that now $v(t)$ is a continuous function
we define
$$
 T(\omega)=\inf\{t\geq 0\mid \|v(t)\|\leq \nu\},\qquad
 T_0=\frac{U_0^\a}{kK}\ln\left(\frac{\|v(0)\|}{\nu}\right).
$$
Taking $\e=\nu\HH_0$ we now see that $T(\omega)\leq T_0$ with probability at least
\begin{equation}\label{eq:bp}
 \left\{
 2\bPhi\left(\frac{\nu\HH_0}{\sigma\sqrt{3k}}\right)
 -\frac{2T_0\sigma\|\dot{\psi}\|}{\delta\sqrt{2\pi}} 
 \varphi\left(\frac{\nu\HH_0}{\sigma\sqrt{3k}}\right)-1
 \right\}^{3k}.
\end{equation}
Let us then take $\omega$ in the set $\max_{0\leq t\leq
T_0}\|H\left(t\right)\|\leq\nu\HH_0$. By our previous computation,
this set has a probability not smaller than the bound
\eqref{eq:bp}. If $T(\omega)>T_0$, then $\|v(t)\|>\nu$ on the
interval $[0,T_0]$, Theorem~\ref{th:1} holds for $T=T_0$, and we obtain
$\|v(t)\|\leq\nu$ for some $t\leq T_0$, obtaining a contradiction.
This concludes the proof. 
\eproof 
\medskip

{\small

}
\end{document}